\documentclass[a4paper,10pt]{article}
\usepackage[top=1in, bottom=1in, left=1in, right=1in]{geometry}
\usepackage{amsfonts}
\usepackage{amssymb}
\usepackage{amsthm}
\usepackage{amsmath,amssymb,latexsym} 
\usepackage{fancyhdr}
\usepackage{graphicx, graphics}

\newtheorem{thm}{Theorem}
\newtheorem{deff}{Definition}
\newtheorem{prop}[thm]{Proposition}
\newtheorem{lem}[thm]{Lemma}
\newtheorem{cor}[thm]{Corollary}
\newtheorem{con}{Conjecture}

\begin{document}

\title{Some divisibility properties of Jacobsthal numbers}
\author{\textit{Volkan Yildiz} }
\date{vo1kan@hotmail.co.uk}
\maketitle

\abstract
In this paper we investigate some divisibility properties of Jacobsthal numbers.
\\\\\\


\pagebreak
\section{Introduction}

Jacobsthal numbers were studied by \cite{H}, and by many others very recently. 
The following observations were made in \cite{H}:
\[
\{J_k\}_{k\ge 0}:= 0, 1, 1, 3, 5, 11, 21, 43, 85, \dots 
\]
It can also be defined as a second order recurrence relation:
\[
J_{k+2} = J_{k+1}+2J_{k}, \;\;\;\; \text{ with }\;\; J_0=0,   \; J_1= 1.
\]
and it's Binet form:
\[
J_k=\frac{2^{k}-(-1)^k}{3}.
\]
There are many subsequences and sibling sequences that one can study among Jacobsthal numbers.
Some are famous and some are not. They all play their parts. Our main approach to these sequences will be looking at their Binet forms:
\[
j_k= 2^k+(-1)^k
\]
Where $j_k$ is known as the Jacobsthal-Lucas sequence and 
it has many interrelationships characteristics which are outlined in \cite{H}. \\\\
Jacobsthal numbers play an important role in substructure of counting the truth or false entires of bracketed formulea connected by the implication operator, see \cite{V1}. I have encountered this interesting sequence while doing calculations with truth tables. Hopefully more to be said in the near future. This paper is consists of a few observations about the divisibility of Jacobsthal numbers.
\begin{deff}
\[
L_{\mu M}=\left\{ \frac{2^{\mu M}+1}{3} : M=2k+1,\; \mu=2l+1,\;\; k,l\in\mathbb{Z}^+ \right\}
\]

\[
L_{\mu N}=\left\{ \frac{2^{\mu N}-1}{3} : N=2k,\; \mu, k\in\mathbb{Z}^+ \right\}
\]
\end{deff}
Working backwards gives:
\begin{equation}\notag
\begin{split}
2^{3(2k+1)}+1&= 2^{3M}+1\\
&= 2^{3M}-2^{2M}+2^M+2^{2M}-2^{M}+1\\
& =(2^M+1)(2^{2M}-2^M+1)\\
\frac{2^{3M}+1}{3}& = 3\Bigg(\frac{2^{M}+1}{3}\Bigg)\Bigg(\frac{2^{2M}-2^{M}+1}{3}\Bigg)\\
J_{3M} & = 3J_{M}a_{M}
\end{split}
\end{equation}
Thus $a_M$ is defined by \cite{S}, A345963. Thus
\[
a_M=\frac{J_{3M}}{3J_{M}}.
\]

Furthermore,
\[
\frac{J_{3M}}{J_M}=3\Big( J_{2M}-J_M+\frac{1}{3} \Big)
\]
and this proves the following proposition.
\[
J_{3M}=J_M\Big(3J_{2M}-3J_{M}+1\Big)
\]
\begin{prop}
\[
\frac{J_{3M}}{J_M}=\Big( 3J_{2M}-3J_M+1  \Big)
\]
\end{prop}\;\\\\
For $k$ is even the following result was mentioned in \cite{H}, and also named as A014551 in \cite{S}
\[
J_{2k}=J_kj_k
\]
It takes only a line to see this:
\[
J_{2k}=\frac{2^{2k}-(-1)^{2k}}{3}=\frac{(2^{k}+1)}{3}(2^{k}-1)=\frac{(2^{k}-1)}{3}(2^{k}+1)
\]
Next observation is even more intersting. Lets consider an odd power of 2:
\begin{prop}
\[
(2^{ab}+1)= (2^a+1) \sum_{i=0}^{b-1} (-2)^{i a}
\]
\end{prop}
\begin{proof}
\begin{equation}\notag
\begin{split}
& (2^a+1)\Big(1-2^a+2^{2a}-2^{3a}+2^{4a}-2^{5a}+\dots -2^{(b-2)a}+2^{(b-1)a}\Big)\\
&=  2^a -2^{2a}+2^{3a}-2^{4a}+\dots -2^{(b-1)a}+2^{ba}\\
&+1-2^a+2^{2a}-2^{3a}+2^{4a}+\dots +2^{(b-1)a}\\
&= (2^{ab}+1).
\end{split}
\end{equation}
\end{proof}
\;\\
Similarly with an even power of 2, we have the following well known result:
\[
(2^{cd}-1)= (2^c-1) \sum_{i=0}^{d-1} 2^{ic}
\]

\section{Divisibility}
In this section we would like to use eliptic divisibility properties of Lucas typed sequences:
\[
n|m \text{ whenever } J_n|J_m.
\]
Sequences satisfying this property are known as divisibility sequences. We also know that sequences satisfying
\[
\gcd(J_n,J_m)=J_{\gcd(m,n)}
\]
are known as strong divisibility sequences; in fact every strong divisibility sequence is a divisibility sequence. 
\[
\gcd(m,n)=m \; \text{ then}  \; m|n.
\]
\begin{prop}
Let $m, n\in \mathbf{Z}^+$ then $n|m$ iff $J_n|J_m$.
\end{prop}

\begin{proof}
 see any elliptic divisibility book.
\end{proof}

\begin{cor}
Let $J_m$ be a Jacobsthal number, and $a_1,\ldots,a_z$ be proper divisors of $m$ then $J_{a_1},\ldots,J_{a_z}$ are Jacobsthal.
\end{cor}

\begin{lem}
Let $q>4$ be composite integer then $J_q$ can be written as product of non-Jacobsthal and Jacobsthal primes.
\end{lem}

\begin{proof}
Let $D$ be the set of divisors of $q$ then let $m=\max(D)$ then 
\[
J_q=J_mJ_k R,\;\;\;\; \text{ where} \;\;\; k=2,...,m-1
\]
Where $R$ is the non-Jacobsthal divisor.\\
Suppose $J_q$ can be written product of Jacobsthal primes only, or non-Jacobsthal primes only. 
The former is not possible, as $q>\lfloor \frac{q}{2}\rfloor =m$ simulteneously $q>k=2$.
By elliptic divisibility the later argument is not possible as $q$ is composite. 
\end{proof}
\begin{cor}
$5=J_4$ is the only Jacobsthal prime with composite index. 
\end{cor}
\begin{con}
There are infinitely many composite Jacobsthal numbers with prime index.
\end{con}

\begin{prop} For $n\geq 0$
\[
 3^n|J_{3^n}
\]
\end{prop}
\begin{proof}
Recall that $\;J_3=3$
\[
J_{\underbrace{3\times \ldots \times 3}_{n\textit{-times}}}=A\times \prod^n J_{3}
\]
for some $A\in \mathbf{Z^+}$
\end{proof}

\begin{cor}
Let $M$ be a factor of $J_{3^k}$ for some $k<n$ then $M\times 3^n$ is a factor of $J_{M\times 3^n}$.
\end{cor}
E.g. 
\[
J_{171}=997718451084563058827048845800982541418449949338283, \;\text{ and} \; 171|J_{171}    
\]
So indices which are in the form of $19\times 3^n$ divides $J_{19\times 3^n}$ for $n>1$. Moreover,
\[
19^n \times 3^{n+2} |J_{(19^n\times 3^{n+2})}, \; \text{ for }  n\in\mathbf{Z^+}.
\]

\begin{prop}
Let $q>4$ and $J_q$ and $q$ be both composite then the maximal divisor of $J_q$ is non-Jacobsthal.
\end{prop}
\begin{proof}
By Lemma 5 we can have:
\[
J_q=J_mJ_kR
\]
Let $a\in\mathbf{Z^+}$
\[
J_a\sim 2^{a-1}, \;\; \text{ for} \;\;a>3, \;\;\ \text{ and} \;\; J_2=J_1=1.
\]
Then
\[
R\sim 2^{{mk-1}-(m-1)-(k-1)}
\]
Thus it is suffice to show that $mk-m-k+1>m-1$. Remember $k=2,...,m-1$ and $m>3$, as
\[
(m-1)(k-2)>1
\]
$R$ is the maximal divisor of $J_q$.
\end{proof}

\begin{prop}
Let $p$ be a prime number greater than 3, then 
\[
J_p\equiv_p 1.
\]
\end{prop}

\begin{proof}
 For $p$ prime, $J_p$ may or may not be a prime number,
but for $p>3$, we know that $J_p$ is not divisible by $p$.
\[
2^p\equiv_p 2 \;\;\Rightarrow\;\; 2^p+1\equiv_p 3 \;\;\Rightarrow\;\; J_p\equiv_p 1.
\]

\end{proof}

\begin{prop}
Let $J_q$ be a composite Jacobsthal, with $J_q\equiv_q M$, and $M\not= 0,1$ then
\[
\exists Q,x\in \mathbf{Z_q^\times} \;\;\; : \;\;\; J_q\equiv_Q 0\;\;\; \wedge \;\;\; xQ\equiv_q 1
\]
\end{prop}

\begin{proof}
By above Lemma, existence of such $Q>q$ can be deduced: 
\[
J_q=J_mJ_k R,\;\;\;\; \text{ where} \;\;\; k=2,...,m
\]
It is clear that $Q$ must be a divisor of $R$\\
\[
\text{with } \; R\sim 2^{(m-1)(k-1)}
\]
We have $\gcd(Q,q)=1$ and by Bezout's Lemma:
\[
xQ+yq=1,\;\;\; \exists x,y\in\mathbf{Z}
\]
Thus
\[
xQ\equiv_q 1 .
\]
\end{proof}
\begin{con}
In above proposition $x=1$.
\end{con}

\pagebreak 

\;\\\\\\\\\\\\\\\\\\\\\\\\\\\\\\\\\\\\\\\\\\\\\\\\\\\\\\\\\\\\\\\\\\\\\\\\

\begin{flushright}
işte yüzünde badem çiçekleri\\
saçlarında gülen toprak ve ilkbahar.\\
sen misin seni sevdiğim o kavga,\\
sen o kavganın güzelliği misin yoksa…\\

A Y\"ucel.
\end{flushright}


\begin{thebibliography}{9}
\bibitem{H}
A F Horadam
\textit{Jacobsthal representation numbers}
The Fibonacci Quarterly, ...


\bibitem{S} 
N. J. A. Slone. 
\textit{OEIS}.
www.OEIS.org


\bibitem{EEW} 
 M Einsiedler,  G Everest and T Ward .\;\\
\textit{Primes in elliptic divisibility sequences,}
LMS, ISSN 1461-1570


\bibitem{div}
Jean-Paul Bézivin, Attila Pethö and Alfred J. van der Poorten\\
A Full Characterisation of Divisibility Sequences\\
AJM,  http://www.jstor.org/stable/2374733
\bibitem{V1} 
V Yildiz \footnote{ Volkan Yildiz is a secondary school teacher, currently living in London. \\ e-mail: vo1kan@hotmail.co.uk. \\ website: https://sites.google.com/site/vo1kanyildizmaths}, \textit{General Combinatorical Structure of Truth Tables of Bracketed Formulae Connected by Implication}.
Arxiv: https://arxiv.org/abs/1205.5595



\end{thebibliography}
\end{document}